\documentclass[a4paper, final, 12pt]{article}
\usepackage{eurosym}
\usepackage{amsmath, amssymb, latexsym, amscd, amsthm,amsfonts,amstext}
\usepackage[mathscr]{eucal}
\usepackage{graphicx}
\usepackage{subfig}
\usepackage{float}
\usepackage{color}
\usepackage{hyperref}
\usepackage[utf8]{inputenc}
\usepackage[english]{babel}

\newtheorem{theorem}{Theorem}[section]

\numberwithin{equation}{section}
\begin{document}

\title{Application of Neural Network Machine Learning to Solution of Black-Scholes Equations}
\author{Mikhail V. Klibanov, Kirill V. Golubnichiy and Andrey V. Nikitin
\and Department of Mathematics and Statistics, \and Department of
Mathematics \and University of North Carolina at Charlotte, Charlotte, USA
\and University of Washington, Seattle, USA \and mklibanv@uncc.edu,
kgolubni@math.washington.edu, \and andrey.nikitin@hotmail.com}
\date{}
\maketitle


\begin{abstract}
This paper presents a novel way to predict options price for one day in
advance, utilizing the method of Quasi-Reversibility for solving the
Black-Scholes equation. The Black-Scholes equation is solved forwards in
time, which is an ill-posed problem. Thus, Tikhonov regularization via the
Quasi-Reversibility Method is applied. This procedure allows to forecast
stock option prices for one trading day ahead of the current one. To enhance
these results, the Neural Network Machine Learning is applied on the second
stage. Real market data are used. Results of Quasi-Reversibility Method and
Machine Learning method are compared in terms of accuracy, precision and
recall.
\end{abstract}

\vspace{0.5em} 




\textbf{\text{Keywords:}}

The Black-Scholes equation, Ill-posed problem, regularization method,
parabolic equation with the reversed time, Machine Learning, neural network.

\section{Introduction}

This paper discusses a new empirical mathematical model for generating more
accurate option trading strategy using initial and boundary conditions for
the underlying stock. The idea was initially proposed in \cite{KlibGol}. The
basis for this idea is the Black-Scholes equation. In mathematical finance,
the Black-Scholes equation is a parabolic partial differential equation that
determines the dynamics of the price of European options \cite{Shreve}.

The time at a given time $t$ will occur is $\tau ,$ 
\begin{equation}
\tau =T-t.  \label{1.0}
\end{equation}%
$f(s)$ be the payoff function of that option
at the maturity time $t=T$ and $s$ is the stock price. Let's assume that the risk-free interest rate equals zero. The function $u\left( s,\tau \right) $\ is the price of that option and the
variable $\tau $\ is the one defined in (\ref{1.0}). Let's assume that this
function $u\left( s,\tau \right) $ satisfies the Black-Scholes equation with
the volatility coefficient ${\sigma}$ \cite[Chapter 7, Theorem 7.7]%
{Bjork}\textbf{:}

\begin{equation}
\begin{split}
& \frac{\partial u(s,\tau )}{\partial \tau }=\frac{{\sigma }^{2}}{2}s^{2}\frac{%
\partial ^{2}u(s,\tau )}{\partial s^{2}}, \\
& u(s,0)=f(s),
\end{split}
\label{1}
\end{equation}%
The payoff function is $f(s)=\max \left(
s-K,0\right) $, where $K$ is the strike price \cite{Bjork} and $s>0.$ 

The option price function is defined by the Black-Scholes formula: 
\begin{equation}
u(s.\tau )=s\Phi (\theta _{+}(s,\tau ))-e^{-r\tau }K\Phi (\theta _{-}(s,\tau
)),
\end{equation}%
%

Based on the It\^o formula, we have:

\begin{equation}
du = (-\frac{\partial u(s, T-t)}{\partial \tau} + \frac{{\sigma}^2}{2} s^2 \frac{%
\partial^2 u(s, T-t)}{\partial s^2})dt + \sigma s \frac{\partial u(s, T-t)}{%
\partial s} dW.
\end{equation}

If equation (\ref{1}) is solved forwards in time to forecast prices of stock
options is an ill-posed inverse problem. For this reason, we used the Method of Quasi- Reversibility (QRM) that is a version of the Tikhonov regularization method. 
Uniqueness, stability and convergence theorems for this method were
formulated in \cite{KlibGol} and \cite{Klibyag}, also, see \cite{SH} for
proofs.

We have four major questions that we raise in this paper:

\begin{enumerate}
\item What is the forecast interval of the options prices?

\item What are the boundary and initial conditions on the interval for the
Black-Scholes equation?

\item What are the values of the volatility coefficient in the future?

\item How to solve the Black-Scholes equation forwards in time $t$? 
\end{enumerate}

The first three questions are addressed in our new mathematical model. We
use the regularization method of \cite{Kl} to address the fourth question.
Theorems about stability and convergence of this method are formulated.
These theorems were proven in \cite{Kl} for a general parabolic equation of
the second order where the main key of this method is based on the method of
Carleman estimates.

Let the function $f$ $\in L_{2}(0,\pi )$ and let $T=const.>0.$ To
demonstrate that our problem is ill-posed, we consider the example of based
on the problem for the heat equation with the reversed time

\begin{equation}
u_{t}+u_{xx}=0,\left( x,t\right) \in \left( 0,\pi \right) \times \left(
0,T\right) ,  \label{2}
\end{equation}%
with Dirichlet boundary conditions 
\begin{equation}
u\left( 0,t\right) =0,u\left( \pi ,t\right) =0,  \label{3}
\end{equation}%
and the initial condition 
\begin{equation}
u\left( x,0\right) =f\left( x\right) .  \label{4}
\end{equation}

The unique solution of this problem is:

\begin{equation*}
u(x,t)=\sum_{n=1}^{\infty }f_{n}\sin {(nx)}e^{n^{2}t}.
\end{equation*}%
Consider%
\begin{equation*}
u_{N}(x,t)=\sum_{n=1}^{N}f_{n}\sin {(nx)}e^{n^{2}t}.
\end{equation*}%
Then 

\begin{equation*}
||u_{N}(x,T)||_{L_{2}(0,\pi
)}^{2}=\sum_{n=1}^{N}f_{n}^{2}e^{2n^{2}T}\approx
f_{N}^{2}e^{2N^{2}T}\rightarrow \infty 
\end{equation*}%
as $N,T\rightarrow \infty .$ Hence,  the problem (\ref{2})-(\ref{4}) is
severely unstable. 

We conclude therefore that to obtain a more or less accurate solution of the
Black-Scholes equation forwards in time, we need to solve it on a short time
interval $\left( 0,T\right) $. To get better accuracy, the regularization
method works only for a short time interval. 

Section 3 presents our mathematical model, the method of Quasi Reversibility
as well as the trading strategy. The Quasi Reversibility Method is based on
the minimization of a Tikhonov-like functional $J_{\beta }(u).$ We do this
using conjugate gradient method. The minimization process was performed by
Hyak Next Generation Supercomputer of the research computing club of
University of Washington. The code was parallelized in order to maximize the
performance on supercomputer clusters.

The historical data for stock options was collected from the Bloomberg
terminal \cite{Bloom} of University of Washington. From this data, we
obtained about 177,000 minimizers.

Due to ill-posedness of the problem the solution is very sensitive to the
noise in the initial data (stock and option prices for the three days
preceding the day of forecast). Given results of the  Quasi Reversibility
Method, we apply on the second stage Machine Learning to reduce the
probability of non-profitable trades caused by wrong option price prognosis
because of the noise in input data.

Section 4 is dedicated to application of binary classification and
regression Neural Network Machine Learning.

Sections 5 and 6 present our results and the summary.

Python with the SciPy and Torch modules were used for implementation of the
method of Quasi-Reversibility and Neural Network Machine Learning (binary
classification and regression).

\section{The new mathematical model and the method of Quasi-Reversibility}

Let's denote $s$ as the stock price, $t$ as the time and $\sigma (t)$ as the
volatility of the option. The historical implied volatility listed on the
market data of \cite{Bloom} is used in our particular case. We assume that $%
\sigma = \sigma(t)$ to avoid other historical data for the volatility. Let's
call $u_{b}(t)$ and $u_{a}(t)$ the bid and ask prices of the options at the
moment of time $t$ and $s_{b}(t)$ and $s_{a}(t)$ the bid and ask prices of
the stock at the moment of time $t.$ It is also known that

\begin{equation}
u_{b}(t)<u_{a}(t)
\end{equation}

and

\begin{equation}
s_{b}(t)<s_{a}(t)
\end{equation}

Let's introduce

\begin{equation}
f_{s}(t)=\frac{s_{a}(t)}{s_{b}(t)}-1
\end{equation}
and 
\begin{equation}
f_{u}(t)=\frac{u_{a}(t)}{u_{b}(t)}-1
\end{equation}

Based on real market data we have observed that usually 
\begin{equation}
0\leq f_{s}(t) \leq 0.003
\end{equation}

and

\begin{equation}
0\leq f_{u}(t) \leq 0.27
\end{equation}

The idea is to approximate the Black-Scholes equation solutions%
\begin{equation}
Lu=u_{t}+\frac{\sigma ^{2}\left( t\right) }{2}s^{2}u_{ss}=0,\left(
s,t\right) \in \left( s_{b}\left( 0\right) ,s_{a}\left( 0\right) \right)
\times \left( 0,2\tau \right) =X_{2\tau },  \label{3.3}
\end{equation}%
with Dirichlet boundary conditions 
\begin{equation}
u\left( s_{b},t\right) =u_{b}\left( t\right) ,u\left( s_{a},t\right)
=u_{a}\left( t\right) ,\text{ }t\in \left[ 0,2\tau \right] ,  \label{3.4}
\end{equation}%
and the initial condition 
\begin{equation}
u\left( s,0\right) =f\left( s\right) ,\text{ \ \ }s\in \left[ s_{b}\left(
0\right) ,s_{a}\left( 0\right) \right] .  \label{3.5}
\end{equation}

Where $L$ is the partial differential operator of the Black-Scholes
equation. Based on Bloomerg terminal we used with End of Day Underlying
Price Last, End of Day Underlying Price Bid, End of Day Underlying Price
Ask, $t$ is time, $\sigma \left( t\right) $ is the volatility of the stock
option. It was used Implied Volatility Using Last Trade Price (IVOL). 

$u\left( s,t\right) $ is the price of the stock option. End of Day Option
Price Last, End of Day Option Price Bid and End of Day Option Price Ask are
the notation that we applied in our algorithm.

\textbf{Problem 1.} \emph{Find the function }$u\in H^{2}\left(
X_{2\tau}\right) $\emph{\ satisfying conditions (\ref{3.3})-(\ref{3.5}).}

This problem considers as ill-posed since we solve equation (\ref{3.3})
forwards in time.

\textbf{Remarks 3.1: }\emph{We increase here the required smoothness of the
solution from }$H^{2,1}\left( X_{2\tau }\right) $\emph{\ to }$H^{2}\left(
X_{2\tau }\right) .$

Our algorithm based on solving the inverse problem for the Black-Scholes
with reversed time equation has five steps:

\textbf{Step 1 (Dimensionless variables).} \bigskip 

We require to make our equation dimensionless. $s_{b} < s_{a}.$ Let's denote 
$s_{b} = s_{b}(0),$ $s_{a} = s_{a}(0).$ Dimensionless variables were applied 
$x,t^{\prime }$ such that 
\begin{equation}
x=\frac{s-s_{b}}{s_{a}-s_{b}}
\end{equation}

\begin{equation}
t^{\prime }=\frac{t}{255}
\end{equation}

and now we can say that $s$ is $x$ and $t$ is $t^{\prime }.$

According to these substitutions, the equation becomes

\begin{equation}
Ru=u_{t}+\sigma ^{2}\left( t\right) A(x)u_{xx},  \label{3.6}
\end{equation}

where 
\begin{equation}
A(x)=\frac{255}{2}\frac{[x(s_{a}-s_{b})+s_{b})]^{2}}{(s_{a}-s_{b})^{2}}
\label{3.7}
\end{equation}

\begin{equation}
X_{2\tau}=\left\{ \left( x,t\right) \in \left( 0,1\right) \times \left(
0,2\tau\right) \right\} .  \label{3.8}
\end{equation}

\begin{equation}
u(x,0)=g(x), x\in(0,1)  \label{3.9}
\end{equation}

\begin{equation}
u(0,t)=u_{b}(t), u(1,t)=u_{a}(t).  \label{3.10}
\end{equation}

And the operator $L$ in (\ref{3.3}) is the operator $R$

\textbf{Step 2 (Interpolation and extrapolation).} \bigskip

Our goal is to forecast option price from 'today' to 'tomorrow' and 'the day
after tomorrow'. We do have 255 trading days annually. For this reason,
let's introduce $\tau >0$ as our unit of time for which we want to make our
prediction the option price. Because we predict option prices having the
information of these prices, as well as of other parameters for 'today',
'yesterday' and 'the day before yesterday', we consider $\tau $ is one
trading day. 'One day' $\tau =1/255.$ 'Today' $t=0.$ 'Tomorrow' $t=\tau .$
'The day after tomorrow' $t=2\tau .$ The variable $s-$ is for interval, i.e $%
s\in \left[ s_{b}\left( 0\right) ,s_{a}\left( 0\right) \right] .$ We applied
the idea associated with interpolation discrete values of functions $%
u_{b}(t),$ $u_{a}(t),$ and $\sigma (t)$ between these three points (the day
before yesterday, yesterday and today) and then extrapolation functions $%
u_{b}(t),$ $u_{a}(t)$ between three points (today, tomorrow and the day
after tomorrow). Where $t=-2\tau $ is "the day before yesterday", $t=-\tau $
is "yesterday" and $t=0$ is "today". We used quadratic polynomials for both
approximation and extrapolation of values of functions. Thus, these three
functions $u_{b}(t),$ $u_{a}(t),$ and $\sigma (t)$ was obtained for a small
future time interval, i.e $(0,2\tau ).$ (\cite{KlibGol}). Where $u_{b}(t),$ $%
u_{a}(t)$ were applied for boundary conditions and $\sigma (t)$ is
coefficient function for our problem. The initial condition was set as $%
u(x,0)=g(x)=x(u_{a}(0)-u_{b}(0))x+u_{b}(0).$ This function is the result of
approximation by linear function due to the fact that the interval between
bid and ask prices is relatively small. The domain was $X_{2\tau
}=\{(x,t):x\in (0,1),t\in (0,2\tau )\}.$ \bigskip 

\textbf{Step 3 (Statement of the Problem).} \bigskip

\textbf{Problem 2.} \emph{Assume that functions }%
\begin{equation}
u_{b}\left( t\right) ,u_{a}\left( t\right) \in H^{2}\left[ 0,2\tau \right]
,\sigma \left( t\right) \in C^{1}\left[ 0,2\tau \right] .  \label{3.11}
\end{equation}%
\emph{Find the solution }$u\in H^{2}\left( X_{2\tau }\right) $\emph{\ of the
following initial boundary value problem:}%
\begin{equation}
Ru=0\text{ in }X_{2\tau },  \label{3.12}
\end{equation}%
\begin{equation}
u\left( 0,t\right) =u_{b}\left( t\right) ,u\left( 1,t\right) =u_{a}\left(
t\right) ,t\in \left( 0,2\tau \right) ,  \label{3.13}
\end{equation}%
\begin{equation}
u\left( x,0\right) =g\left( x\right) ,x\in \left( 0,1\right) ,  \label{3.14}
\end{equation}%
\emph{where the partial differential operator }$R$\emph{\ is defined in (\ref%
{3.6}), the function }$A\left( x\right) $\emph{\ is defined in (\ref{3.7}),
the initial condition }$g\left( x\right) $\emph{\ is defined in (\ref{3.9}),
and the domain }$X_{2\tau }$\emph{\ is defined in (\ref{3.8}).}

\begin{theorem}
The following problem (\ref{3.6})-(\ref{3.9}) has one solution $u \in
H^{2,1}(X_{2\tau}).$
\end{theorem}

The proof of this theorem is \cite{Klibyag}. \bigskip

\textbf{Step 4 (Numerical method of solving the problem. Regularization).}
\bigskip

Due to the ill-posedness of the problem, we can not say about existence of
the solution. Thus, it was applied the regularization method:

Let's consider function $F(x,t)=x(u_{a}(t)-u_{b}(t))+u_{b}(t),$ $(x,t)\in
X_{2\tau}.$ This function $F \in H^{2}(X_{2\tau}).$ It follows from (\ref%
{3.9}) and (\ref{3.10}) that

\begin{equation}
F\left( x,0\right) =g\left( x\right) ,  \label{4000}
\end{equation}

\begin{equation}
F(0,t ) =u_{b}\left( t\right) ,F(1,t) =u_{a}\left( t\right).  \label{3000}
\end{equation}

We used an unbounded differential operator $R: H^{2,1}(X_{2\tau})
\rightarrow L^{2}(X_{2\tau}),$ where $H^{2,1}(X_{2\tau})$ is a dense linear
set in the space $L^{2}(X_{2\tau}).$ Where 
\begin{equation}
Ru=u_{t}+\sigma ^{2}\left( t\right) A(x)u_{xx}
\end{equation}

Let's introduce Tikhonov-like functional as:

\begin{equation}
J_{\beta }\left( u\right) =\int_{X_{2\tau }}\left( Ru\right) ^{2}dsdt+\beta
\left\Vert u\right\Vert _{H^{2}\left( X_{2\tau }\right) }^{2},  \label{3.15}
\end{equation}

where $\beta \in \left( 0,1\right) $ is the parameter of regularization. To
solve the problem, we minimized the functional $J_{\beta }\left( u\right) $
on the set

\begin{equation}
V=\left\{ u\in H^{2}\left( X_{2\tau}\right) :u\left( 0,t\right) =u_{b}\left(
t\right) ,u\left( 1,t\right) =u_{a}\left( t\right) ,u\left( x,0\right)
=g\left( x\right) \right\} .  \label{3.16}
\end{equation}

\textbf{Step 5 (Minimization Problem).} \bigskip

\textbf{Minimization Problem 1}. $J_{\beta }:H^{2}\left( X_{2\tau}\right)
\rightarrow \mathbb{R}$\emph{\ is the regularization Tikhonov functional.}

We have used the converting of our partial derivatives from (\ref{3.15})
into finite differences. A finite difference grid was applied to cover the
domain $X_{2\tau}.$ The minimization process was to differentiate our
functional $J_{\beta}(u)$ with respect to the values of the function $u(x,t)$
at each grid points via conjugate gradient method. The point $u=0$ was used
for the starting point. Based on computational study with simulated data we
have realized that the optimal value of the regularization parameter would
be $\beta=0.01.$

Minimization Problem 1 is a QRM for Problem 2. This is an version of the QRM
for problem (\ref{3.12})-(\ref{3.14}). In section 4 we discuss the theory of
this specific version of the QRM. In particular, Theorem 4.2 of section 4
presents uniqueness of the solution $u\in H^{2}\left( X_{2\tau}\right) $ of
Problem 2 and implies an estimate of the stability of this solution with
respect to the noise in the data. Theorem 4.3 of section 4 \ shows existence
and uniqueness of the minimizer $u_{\beta }\in H^{2,1}\left(
X_{2\tau}\right) $ of the functional $J_{\beta }\left( u\right) $ on the set 
$V$ defined in (\ref{3.16}). We call such a minimizer \textquotedblleft
regularized solution" \cite{T}. Theorem 4.4 estimates convergence rate of
regularized solutions to the exact solution of Problem 2 with the noiseless
data. Such estimates depend on the noise level in the data. All proof of
these theorems are presented in \cite{SH}.

\section{Analysis}

This section is devoted to convergence analysis for Problem 2 of subsection
3.2. This problem is the initial boundary value problem for parabolic
equation (\ref{3.12}) with the reversed time. The QRM and convergence
analysis for this problem for a more general parabolic operator in $\mathbb{R%
}^{n}$ with arbitrary variable coefficients was proposed in \cite{Kl}. Then
theorems were presented in \cite{KlibGol}. However a stability estimate was
not a part of \cite{KlibGol}, such an estimate was proven in \cite{Kl}. The
same is true for the convergence theorems of QRM\ in \cite{Kl,KlibGol}. The
smallness assumption was lifted in \cite{Klibyag} via a new Carleman
estimate. Results of \cite{Klibyag} for a 1-D case were significantly
modified in this section. Our computations below on a small time interval $%
\left( 0,2\tau\right) =\left( 0,0.00784\right) $ (see \cite{Kl,KlibGol}, 
\cite[Theorem 1 of section 2 in Chapter 4]{LavR} might result in the
requirement of even a smaller length of that interval.

\subsection{Problem statement}

\label{sec:4.1}

Let's consider a number $T>0$ and introduce $Q_{T}$ as: 
\begin{equation*}
Q_{T}=\left\{ \left( x,t\right) \in \left( 0,1\right) \times \left(
0,T\right) \right\} .
\end{equation*}%
Consider two numbers $b_{0},b_{1}>0$ and $b_{0}<b_{1}.$ Let the function $%
b\left( x,t\right) \in C^{1}\left( \overline{Q}_{T_{}}\right) $ satisfies:%
\begin{equation}
\text{ }\left\Vert b\right\Vert _{C^{1}\left( \overline{Q}_{T_{}}\right)
}\leq b_{1},\text{ }b\left( x,t\right) \geq b_{0}\text{ in }Q_{T_{}}.
\label{7.1}
\end{equation}%
We also have functions $\psi _{0}\left( t\right) ,\psi _{1}\left( t\right)
\in H^{2}\left( 0,T_{}\right) .$ In the above case of subsection 3.2, 
\begin{equation*}
T=2\tau,b\left( x,t\right) =\sigma ^{2}\left( t\right) A(x),\psi _{0}\left(
t\right) =u_{b}\left( t\right) ,\psi _{1}\left( t\right) =u_{a}\left(
t\right) .
\end{equation*}%
We now formulate Problem 3, which is a slight generalization of Problem 2.

\textbf{Problem 3}. \emph{Find a solution }$v\in H^{2}\left( Q_{T_{}}\right) 
$\emph{\ of the following (IBVP):}%
\begin{equation}
Nv=v_{t}+b\left( x,t\right) v_{xx}=0\text{ in }Q_{T_{}},  \label{7.3}
\end{equation}%
\begin{equation}
v\left( 0,t\right) =\psi _{0}\left( t\right) ,v\left( 1,t\right) =\psi
_{1}\left( t\right) ,\text{ }t\in \left( 0,T_{}\right) ,  \label{7.4}
\end{equation}%
\begin{equation}
v\left( x,0\right) =z\left( x\right) =\psi _{0}\left( 0\right) \left(
1-x\right) +\psi _{1}\left( 0\right) x,\text{ }x\in \left( 0,1\right) .
\label{7.5}
\end{equation}

\textbf{Remark 4.1.} \emph{Because Problem 2 is less general than Problem 3,
then this analysis of converegence for Problem 3 also works for Problem 2.}

We use the linear function for $v\left( x,0\right) $ in (\ref{7.5}) is to
simplify the initial condition in (\ref{3.14}). Now problem 3 is an IBVP\
for the parabolic equation (\ref{7.3}) with the reversed time. For this
reason, the problem can be considered as ill-posed. Assume that the boundary
with a noise of the level $\nu >0$ in (\ref{7.4}) are in place. Here $\nu $
is a sufficiently small number, i.e.%
\begin{equation}
\left\Vert \psi _{0}-\psi _{0}^{\ast }\right\Vert _{H^{1}\left(
0,T_{}\right) }<\nu ,\left\Vert \psi _{1}-\psi _{1}^{\ast }\right\Vert
_{H^{1}\left( 0,T_{}\right) }<\nu ,  \label{7.6}
\end{equation}%
where functions $\psi _{0}^{\ast },\psi _{1}^{\ast }\in H^{2}\left(
0,T_{}\right) $ are \textquotedblleft ideal" noiseless data. we assume that
there exists an exact solution $v^{\ast }\in H^{2}\left( Q_{T_{}}\right) $
of problem (\ref{7.3})-(\ref{7.5}) with these noiseless data (based on on
the theory of Ill-Posed problems). Below we present estimates how this noise
affects the accuracy of the solution of Problem 3 and also discuss the
convergence rate of numerical solutions obtained by QRM to the exact one as $%
\nu \rightarrow 0.$

Let's introduce the version of functional (\ref{3.15}):%
\begin{equation}
I_{\beta }\left( v\right) =\int_{Q_{T_{}}}\left( Nv\right) ^{2}dxdt+\beta
\left\Vert v\right\Vert _{H^{2}\left( Q_{T_{}}\right) }^{2}.  \label{7.7}
\end{equation}%
We also have the set $W\subset H^{2}\left( Q_{T_{}}\right) ,$%
\begin{equation}
W=\left\{ v\in H^{2}\left( Q_{T_{}}\right) :v\left( 0,t\right) =\psi
_{0}\left( t\right) ,v\left( 1,t\right) =\psi _{1}\left( t\right) ,v\left(
x,0\right) =z\left( x\right) \right\} .  \label{7.8}
\end{equation}%
The solution of Problem 3 is approximate solution by solving the following
problem:

\textbf{Minimization Problem 2}. \emph{Minimize the functional }$I_{\beta
}\left( v\right) $\emph{\ on the set }$W$\emph{\ given in (\ref{7.8}).}

Minimization Problem 2 is QRM for Problem 3.

\subsection{Theorems}

This subsection presents four theorems for Problem 3. All proofs might be
found in \cite{SH}. First, let's introduce the Carleman Weight Function $%
\phi _{\alpha }\left( t\right) $ with $\alpha >2$ for the operator $\partial
_{t}+b\left( x,t\right) \partial _{x}^{2}$ as:%
\begin{equation}
\phi _{\alpha }\left( t\right) =e^{\left( T_{}+1-t\right) ^{\alpha }},\text{ 
}t\in \left( 0,T_{}\right) .  \label{7.9}
\end{equation}%
As a result, the function $\phi _{\alpha }\left( t\right) $ is decreasing on 
$\left[ 0,T_{}\right] $, $\phi _{\alpha }^{\prime }\left( t\right) <0,$%
\begin{equation}
\max_{\left[ 0,T_{}\right] }\phi _{\alpha }\left( t\right) =\psi _{\alpha
}\left( 0\right) =e^{\left( T_{}+1\right) ^{\alpha }},\text{ }\min_{\left[
0,T_{}\right] }\phi _{\alpha }\left( t\right) =\phi _{\alpha }\left(
T_{}\right) =e.  \label{7.10}
\end{equation}%
Denote%
\begin{equation}
H_{0}^{2}\left( Q_{T_{}}\right) =\left\{ u\in H^{2}\left( Q_{T_{}}\right)
:u\left( 0,t\right) =u\left( 1,t\right) =0\right\} .  \label{7.99}
\end{equation}%
\begin{equation}
H_{0,0}^{2}\left( Q_{T_{}}\right) =\left\{ u\in H_{0}^{2}\left(
Q_{T_{}}\right) :u\left( x,0\right) =0\right\} .  \label{7.100}
\end{equation}

\textbf{Theorem 4.1} (Carleman estimate). \emph{Let the coefficient }$%
b\left( x,t\right) $\emph{\ of the operator }$N$\emph{\ satisfies conditions
(\ref{7.1}). Then there exist a sufficiently large number }$\alpha
_{0}=\alpha _{0}\left( T_{},b_{0},b_{1}\right) >2$\emph{\ and a constant }$%
C=C\left( T_{},b_{0},b_{1}\right) >0,$\emph{\ both depending only on listed
parameters, such that the following Carleman estimate holds for the operator 
}$N:$\emph{\ }%
\begin{equation*}
\int_{Q_{T_{}}}\left( Nu\right) ^{2}\phi _{\alpha }^{2}dxdt\geq C\sqrt{%
\alpha }\int_{Q_{T_{}}}u_{x}^{2}\psi _{\alpha }^{2}dxdt+C\alpha
^{2}\int_{Q_{T_{}}}u^{2}\phi _{\alpha }^{2}dxdt
\end{equation*}%
\begin{equation}
-C\sqrt{\alpha }\left\Vert u\right\Vert _{H^{2}\left( Q_{T_{}}\right)
}^{2}-C\lambda \left( T_{}+1\right) ^{\alpha }e^{2\left( T_{}+1\right)
^{\alpha }}\left\Vert u\left( x,0\right) \right\Vert _{L_{2}\left(
0,1\right) }^{2},  \label{7.11}
\end{equation}%
\begin{equation*}
\forall \alpha \geq \alpha _{0},\forall u\in H_{0}^{2}\left( Q_{T_{}}\right)
.
\end{equation*}

Carleman estimate (\ref{7.11}) is the MAIN TOOL to proofs of Theorems 4.2,
4.4.

\textbf{Theorem 4.2} (H\"{o}lder stability estimate for Problem 3 and
uniqueness). \emph{Let the coefficient }$b\left( x,t\right) $\emph{\ of the
operator }$N$\emph{\ satisfies conditions (\ref{7.1}). Let's assume that the
functions }$v\in H^{2}\left( Q_{T_{}}\right) $\emph{\ and }$v^{\ast }\in
H^{2}\left( Q_{T_{}}\right) $\emph{\ are solutions of Problem 3 with the
vectors of data }$\left( \psi _{0}\left( t\right) ,\psi _{1}\left( t\right)
\right) $\emph{\ and }$\left( \psi _{0}^{\ast }\left( t\right) ,\psi
_{1}^{\ast }\left( t\right) \right) $\emph{\ respectively, where }$\psi
_{0},\psi _{1},\psi _{0}^{\ast },\psi _{1}^{\ast }\in H^{2}\left(
0,T_{}\right) .$\emph{\ Assume also that error estimates (\ref{7.6}) of the
boundary data is in place. Choose an arbitrary number }$\epsilon \in \left(
0,T_{}\right) $\emph{. Denote }%
\begin{equation}
\lambda =\lambda \left( T_{},\epsilon \right) =\frac{\ln \left(
T_{}+1-\epsilon \right) }{\ln \left( T_{}+1\right) }\in \left( 0,1\right) .
\label{7.110}
\end{equation}%
\emph{Then there exists a sufficiently small number }$\nu _{0}=\nu
_{0}\left( T_{1},b_{0},b_{1}\right) \in \left( 0,1\right) $\emph{\ and a
constant }$C_{1}=C_{1}\left( T_{},b_{0},b_{1},\epsilon \right) >0,$ \emph{%
both depending only on listed parameters, such that\ the following stability
estimate holds for all }$\nu \in \left( 0,\nu _{0}\right) :$%
\begin{equation}
\left\Vert v_{x}-v_{x}^{\ast }\right\Vert _{L_{2}\left( Q_{T_{}-\epsilon
}\right) }+\left\Vert v-v^{\ast }\right\Vert _{L_{2}\left( Q_{T_{}-\epsilon
}\right) }\leq  \label{7.12}
\end{equation}%
\begin{equation*}
\leq C_{1}\left( 1+\left\Vert v-v^{\ast }\right\Vert _{H^{2}\left(
Q_{T_{}}\right) }\right) \exp \left[ -\left( \ln \nu ^{-1/2}\right)
^{\lambda }\right] .
\end{equation*}

Below $C=C\left( T_{},b_{0},b_{1}\right) >0$ and $C_{1}=C_{1}\left(
T_{},a_{0},b_{1}\right) >0$ denote different constants depending only on
listed parameters.

\textbf{Corollary 4.1} (uniqueness). \emph{Let the coefficient }$b\left(
x,t\right) $\emph{\ of the operator }$N$\emph{\ satisfies conditions (\ref%
{7.1}). Then Problem 3 has at most one solution (uniqueness)}.

\textbf{Proof}.\emph{\ }If $\nu =0,$\ then (\ref{7.12}) implies that $%
v\left( x,t\right) =v^{\ast }\left( x,t\right) $\ in $Q_{T_{}-\epsilon }.$
Since $\epsilon \in \left( 0,T_{}\right) $\ is an arbitrary number, then $%
v\left( x,t\right) \equiv v^{\ast }\left( x,t\right) $\ in $Q_{T_{}}.$ $%
\square $

\textbf{Theorem 4.3 }(existence and uniqueness of the minimizer)\textbf{. }%
\emph{Let functions }$\psi _{0}\left( t\right) ,\psi _{1}\left( t\right) \in
H^{2}\left( 0,T_{}\right) .$\emph{\ Let }$W$\emph{\ be the set defined in (%
\ref{7.8}). Then there exists unique minimizer }$v_{\min }\in W$\emph{\ of
functional (\ref{7.7}) and }%
\begin{equation}
\left\Vert v_{\min }\right\Vert _{H^{2}\left( Q_{T_{}}\right) }\leq \frac{C}{%
\sqrt{\beta }}\left( \left\Vert \psi _{0}\right\Vert _{H^{2}\left(
0,T_{}\right) }+\left\Vert \psi _{1}\right\Vert _{H^{2}\left( 0,T_{}\right)
}\right) .  \label{7.13}
\end{equation}

In the theory of Ill-Posed Problems, this minimizer $v_{\min }$ is called
\textquotedblleft regularized solution" of Problem 3 \cite{T}. According to
the theory of Ill-Posed problems, it is important to establish convergence
rate of regularized solutions to the exact one $v^{\ast }.$ In doing so, one
should always choose a dependence of the regularization parameter $\beta $
on the noise level $\nu ,$ i.e. $\beta =\beta \left( \nu \right) \in \left(
0,1\right) $ \cite{T}.

\textbf{Theorem 4.4} (convergence rate of regularized solutions). \emph{Let }%
$v^{\ast }\in H^{2}\left( Q_{T_{}}\right) $\emph{\ be the solution of
Problem 3 with the noiseless data }$\left( \psi _{0}^{\ast }\left( t\right)
,\psi _{1}^{\ast }\left( t\right) \right) .$\emph{\ Let functions }$\psi
_{0},\psi _{1},\psi _{0}^{\ast },\psi _{1}^{\ast }\in H^{2}\left(
0,T_{}\right) .$\emph{\ Let }$v_{\min }\in W$\emph{\ be the unique minimizer
of functional (\ref{7.7}) on the set }$W$\emph{. Assume that error estimates
(\ref{7.6}) hold. Choose an arbitrary number }$\epsilon \in \left(
0,T_{}\right) $\emph{. Let }$\lambda =\lambda \left( T_{},\epsilon \right)
\in \left( 0,1\right) $ \emph{be the number defined in (\ref{7.110}) and let}%
\begin{equation}
\beta =\beta \left( \nu \right) =\nu ^{2},  \label{7.140}
\end{equation}%
\emph{\ Then there exists a sufficiently small number }$\nu _{0}=\nu
_{0}\left( T_{},b_{0},b_{1}\right) \in \left( 0,1\right) $\emph{\ depending
only on listed parameters such that the following convergence rate of
regularized solutions }$v_{\min }$ \emph{holds for all }$\nu \in \left(
0,\nu _{0}\right) :$\emph{\ }%
\begin{equation}
\left\Vert \partial _{x}v_{\min }-\partial _{x}v^{\ast }\right\Vert
_{L_{2}\left( Q_{T_{}-\epsilon }\right) }+\left\Vert v_{\min }-v^{\ast
}\right\Vert _{L_{2}\left( Q_{T_{}-\epsilon }\right) }  \label{7.14}
\end{equation}%
\begin{equation*}
\leq C_{1}\left( 1+\left\Vert v^{\ast }\right\Vert _{H^{2}\left(
Q_{T_{}}\right) }+\left\Vert \psi _{0}^{\ast }\right\Vert _{H^{2}\left(
0,T_{}\right) }+\left\Vert \psi _{1}^{\ast }\right\Vert _{H^{2}\left(
0,T_{}\right) }\right) \exp \left[ -\left( \ln \nu ^{-1/2}\right) ^{\lambda }%
\right] .
\end{equation*}



\subsection{Trading Strategy:}

We use minimizers obtained from the method of Quasi-Reversibility to build a
strategy for trading options. Let's define

\begin{equation}
REAL (0)= \frac{u_{a}(0)+u_{b}(0)}{2}
\end{equation}

\begin{equation}
REAL (\tau )= \frac{u_{a}(\tau)+u_{b}(\tau)}{2}
\end{equation}

\begin{equation}
EST (\tau)= u_{\beta}(1/2, k\tau)
\end{equation}

or if it was not applied dimensionless

\begin{equation}
EST (\tau)= u_{\beta}(\frac{s_{a}+s_{b}}{2}, k\tau)
\end{equation}

where $k=1$

Here $EST (\tau)$ means minimizer.

Let's buy an option if the following holds

\begin{equation}
EST(\tau) \geq REAL (0)
\end{equation}
The predicted outcome of option trade is Positive if

\begin{equation}
EST(\tau) \geq REAL (0)
\end{equation}
\textbf{Definition 1.}

It is True Positive if 
\begin{equation}
EST(\tau) \geq REAL (0)
\end{equation}
and 
\begin{equation}
REAL(\tau) \geq REAL (0)
\end{equation}
\textbf{Definition 2.}

It is True Negative if

\begin{equation}
EST(\tau)<REAL (0)
\end{equation}
and 
\begin{equation}
REAL(\tau)< REAL (0)
\end{equation}
\newline

\textbf{Definition 3.}

It is False Positive if 
\begin{equation}
EST(\tau) \geq REAL(0)
\end{equation}
and

\begin{equation}
REAL(\tau) < REAL (0)
\end{equation}

\textbf{Definition 4.}

It is False Negative if 
\begin{equation}
EST(\tau)<REAL(0)
\end{equation}
and

\begin{equation}
REAL(\tau) \geq REAL (0)
\end{equation}
\newline
\newline
The accuracy of trading strategy is defined as 
\begin{equation}
Accuracy =\frac{TP+TN}{\sum options}
\end{equation}
where $TP$ is a summation of True Positive and $TN$ is a summation of True
Negative and $\sum options$ is a summation of options in data set.

The precision of trading strategy is defined as 
\begin{equation}
Precision = \frac{TP}{TP+FP}
\end{equation}
where $FP$ is a summation of False Positive.

The recall of trading strategy is defined as 
\begin{equation}
Recall = \frac{TP}{TP+FN}
\end{equation}
where $FN$ is a summation of False Negative.

The average relative error of trading strategy is defined as 
\begin{equation}
Error = \frac{1}{N}\sum |\frac{EST(\tau)-REAL(\tau)}{REAL(\tau)}|
\end{equation}

\newpage

\section{Application of Neural Network Machine Learning}

The Black-Scholes equation gives fair value of options in perfect market.
However, real options prices contain some level of noise. We try to filter
mispredictions (i.e. where minimizers result in False Positive or False
Negative) caused by input noise using Machine Learning to improve accuracy,
precision and recall of the trading strategy. We built a neural network with
13 element input vector and 3 fully connected hidden layers. (See Fig 1).
Input vector consists of minimizers (for $t=\tau, 2\tau$) obtained from the
method of Quasi-Reversibility, stock ask and bid price (for $t=0$), option
ask and bid price and volatility (for $t=-2\tau, -\tau,0$).


\begin{figure}[h!]
\begin{center}
\includegraphics[width =1\textwidth]{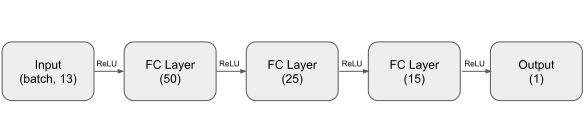}
\end{center}
\par
\caption{}
\end{figure}

All vectors and labels are split into three parts: training, validation and
test sets. The training set is used for weight learning. Validation set is
used for tuning of the neural network hyper-parameters. Test set is for
generating the outcomes of trading strategy.

We collected historical option and stock prices along with implied
volatility on companies consisting of Russel 2000 index \cite{Russ}.

\bigskip

\textbf{Table 1.} \vspace{0.5cm
}

\begin{tabular}{|l|c|l|}
\hline
Set & Dates & Number of options \\ \hline
Training & 2016/09/14-2018/05/31 & 132,912 \\ \hline
Validation & 2018/06/01-2018/06/29 & 13,401 \\ \hline
Test & 2018/07/02-2018/08/17 & 23,549 \\ \hline
\end{tabular}

\bigskip

We compared the profitability of the trading strategy based on the original
minimizer set with the profitability of the output of Machine Learning.

\bigskip

\subsection{Machine Learning Input Vector Normalization}

\begin{equation}
\mu = \frac{u_{a}(0)+u_{a}(-\tau)+u_{a}(-2\tau)+u_{b}(0)+u_{b}(-%
\tau)+u_{b}(-2\tau)}{6}
\end{equation}

\begin{equation}
op_{n}= \frac{u(t)-\mu}{\sigma}
\end{equation}

\begin{equation}
s_{n}= \frac{(s-st)-\mu}{\sigma}
\end{equation}
where $op_{n}$ is a normalized option price, $s_{n}$ is a normalized stock
price normalization, $s$ is the stock, $st$ is the strike and $\sigma$ is
the standard deviation.

\subsection{Binary classification}

\bigskip

Supervised Machine Learning has been applied to the neural network for the
Cross Entropy Loss function with regularization: 
\begin{equation}
L(\theta)=\frac{1}{m}\sum_{i=1}^{m}[-y^{(i)}\log(h_{\theta}(x^{(i)})
-(1-y^{(i)})\log (1-h_{\theta}(x^{(i)}))] + \frac{\lambda}{2m}
\sum_{j=1}^{n}\theta_{j}^{2}
\end{equation}

Where $\theta$ are weights which are optimized by minimizing the loss
function using the method of gradient descent. $\lambda$ is a parameter of
regularization. $x^{(i)}$ is our normalized 13 - dimensional vectors. $%
h_{\theta}$ is output of the neural network. $m$ is the number of vectors in
the training set. $y^{(i)}$ is our labels (the ground truth). The trading
strategy is defined by 
\begin{equation}
H_{c}= 
\begin{cases}
1, & \text{if } h_{\theta}>c \\ 
0, & \text{otherwise}.%
\end{cases}%
\end{equation}
where $c$ is the threshold obtained by maximizing accuracy on validation
set. The labels are set to 1 for profitable trades and 0 otherwise.

\subsection{Regression model}

Similarly, instead of using binary classification, we can use the same ML
architecture to predict the option price for tomorrow ($\tau$). We have the
same input features as the classification neural network. Regression
learning uses mean squared error as the loss function:

\begin{equation}
L(\bar h_{\theta}, \bar y) = \frac{1}{m} \sum_{n} (\bar h_{\theta} - \bar
y_n)^2
\end{equation}

Where $m$ is the size of the data set, $\bar h_{\theta}$ is the predicted
value and $\bar y$ is the real value ($Real (\tau)$).


\section{Results}

The following graph shows the accuracy of the results on validation set. We
use it to determine the optimal value of hyper-parameter $c$ (the threshold
value of binary classification).

\clearpage

\begin{figure}[h]
\par
\begin{center}
\includegraphics[width =1.0\textwidth]{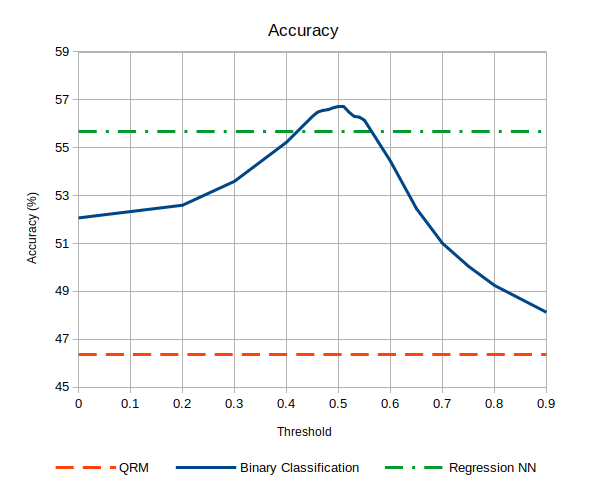}
\end{center}
\par
\caption{Accuracy. Threshold dependency}
\end{figure}

\textbf{Observation 1.}

\bigskip

The accuracy was improved by both Machine Learning methods compared to the
method of Quasi-Reversibility. Based on this graph we set $c=0.5.$

\bigskip

The next graph presents Recall and Precision diagram built on validation
data set.

\clearpage

\begin{figure}[h]
\par
\begin{center}
\includegraphics[width =1.0\textwidth]{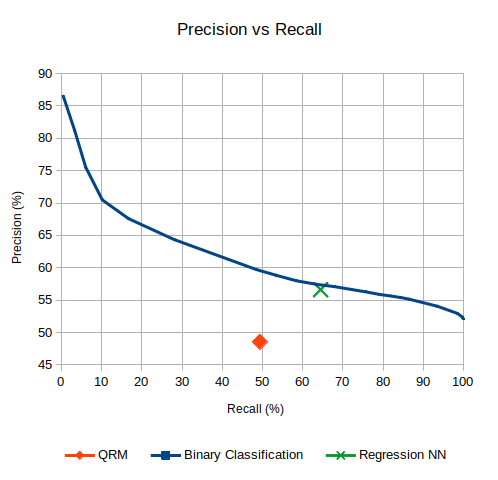}
\end{center}
\par
\caption{Precision vs Recall. Threshold dependency. Here \textbf{X} indicates the position of Recall and Precision produced by Regression Neural Network and $\Diamond$ indicates the position of Recall and Precision produced by QRM.}
\end{figure}


\textbf{Observation 2.}

\bigskip

Binary Classification and Regression produced similar results that improved
both precision and recall compared to the method of Quasi-Reversibility.

\newpage

Further we divided our test data into bins where (horizontal axis, see
Figure ~\ref{fig:Binary Classification, Regression NN and the method of
Quasi-Reversibility.}) each bin is determined by $\frac{s-st}{s}$ with step
size $0.1.$ for each bin we calculated precision (see Figure ~\ref%
{fig:Binary Classification, Regression NN and the method of
Quasi-Reversibility.}).


\begin{figure}[h!]
\par
\begin{center}
{\includegraphics[width = 1.0\textwidth]{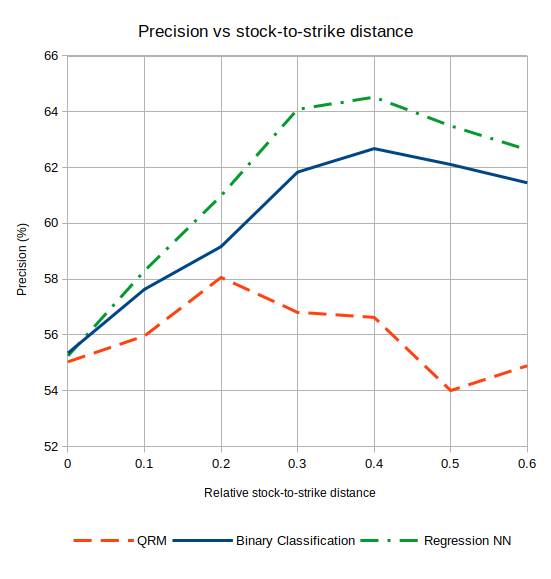}} 
\end{center}
\caption{Binary Classification, Regression NN and the method of
Quasi-Reversibility. }
\label{fig:Binary Classification, Regression NN and the method of
Quasi-Reversibility.}
\end{figure}


\textbf{Observation 3. }

\bigskip

To our surprise, when stock price was close to the strike price Machine
Learning and the method of Quasi-Reversibility give similar precision (bin $%
0 $). With stock price diverging from the strike price Machine Learning
produced better precision.

\clearpage
\newpage

The following tables summarize the accuracy, precision and recall for all
methods on test data.

\bigskip \vspace{0.5cm
} \textbf{Table 2. Final results on Test Data.}%
\vspace{0.5cm
}

\begin{tabular}{|l|l|l|l|l|}
\hline
Method & Accuracy & Precision & Recall & Error \\ \hline
QRM & 49.77\% & 55.77\% & 52.43\% & 12 \% \\ \hline
Binary Classification & 56.36\% & 59.56\% & 70.22\% & NA \\ \hline
Regression NN & 55.42\% & 60.32\% & 61.29\% & NA \\ \hline
\end{tabular}

\bigskip

\bigskip

\textbf{Table 3. Percentages of options with profits/losses for three
different methods. } \vspace{0.5cm
}

\begin{tabular}{|l|c|l|}
\hline
Method & Profitable options & Options with loss \\ \hline
QRM & 55.77\% & 44.23 \% \\ \hline
Binary Classification & 59.56\% & 40.44\% \\ \hline
Regression NN & 60.32\% & 39.68\% \\ \hline
\end{tabular}


\section{Summary}

To predict prices of stock options, we used two empirical mathematical
models for Black-Scholes equation. The results achieved by
solving the equation forwards in time (as an ill-posed problem) and applying
Supervised Machine Learning (Binary Classification and Regression Neural
Network, and using these methods with the real market data, show that this
methodology produce promising results, potential applications within
real-world trading and investment strategies.

The comparison of our methods resulted in the following two conclusions:

\begin{enumerate}
\item The predictions of the method of Quasi-Reversibility ended up being
profitable for $55.77\%$ of the options. Compare this to a $59.56\%$
profitability rate for the Binary Classification method, and a $60.32\%$
profitability rate for Regression Neural Network, used on the same data set
and with the same trading strategy.

\item As shown by figures in section 5, option price forecasting using
Machine Learning gives us significant accuracy and profit improvements over
the method of Quasi-Reversibility. However, when stock price is close to the
strike price both models give similar results.
\end{enumerate}

The authors hypothesize that options traders can generate significant
profits using trading strategies reliant on predictions generated with these
methods.


\newpage



\begin{thebibliography}{99}
\bibitem{Bakush} 

\textsc{M.V. Klibanov A.B. Bakushinskii and N.A. Koshev}, \emph{Carleman weight
functions for a globally convergent numerical method for ill-posed Cauchy
problems for some quasi- linear PDEs}, Nonlinear Analysis: Real World
Applications, 34:201–224, 2017.

\bibitem{Bloom} 

\url{https://bloomberg.com.}

\bibitem{Klibyag} 

\textsc{M. V. Klibanov and A. G. Yagola}, \emph{Convergent numerical methods for
parabolic equations with reversed time via a new Carleman estimate}, Inverse
Problems, 35, 2019.

\bibitem{Kl} 

\textsc{M.V. Klibanov,} \emph{Carleman estimates for the
regularization of ill-posed Cauchy problems}, Appl. Numer. Math., 94, 46-74,
2015.

\bibitem{KlibGol} 

\textsc{M.V. Klibanov A.V. Kuzhuget and K.V. Golubnichiy}, \emph{An ill-posed problem
for the Black-Scholes equation for a profitable forecast of prices of stock
options on real market data}, Inverse Problems, 32(1), 2016.

\bibitem{Russ} 

\url{https://money.cnn.com/data/markets/russell.}

\bibitem{Shreve} 

\textsc{S. E. Shreve}, \emph{Stochastic Calculus for Finance II. Continuous - Time
Models}, Springer, 2003

\bibitem{LavR} 

\textsc{M.M. Lavrent'ev, V.G. Romanov and S.P. Shishatskii}, \emph{Ill-Posed Problems
of Mathematical Physics and Analysis}, Providence, RI: American Mathematical
Society, 1986.

\bibitem{T} \textsc{A.~N. Tikhonov, A.~V. Goncharsky, V.~V. Stepanov and
A.~G. Yagola}, \emph{Numerical Methods for the Solution of Ill-Posed Problems%
}, Kluwer Academic Publishers Group, Dordrecht, 1995.

\bibitem{SH} \textsc{M.V. Klibanov, A. A. Shananin, K. V. Golubnichiy and S.
M. Kravchenko}, \emph{Forecasting Stock Options Prices via the Solution of
an Ill-Posed Problem for the Black-Scholes Equation}, arXiv preprint
arXiv:2202.07174.

\bibitem{Bjork} \textsc{T. Bjork}, \emph{Arbitrage Theory in Continuous Time}%
, Oxford University Press, 1999.



\end{thebibliography}
\end{document}